\documentclass[12pt]{article}
\usepackage{amssymb, amsmath}
\newcommand{\eq}{\begin{equation}}
\newcommand{\en}{\end{equation}}
\newcommand{\te}{\rightarrow}
\newcommand{\re}[1]{\mbox{(\ref{#1})}}

\newcommand{\nutil}{\tilde{\nu}}
\newcommand{\prob}{\mathbb P}
\newcommand{\N}{\mathbb N}
\newcommand{\giv}{\,|\,}
\newcommand{\B}{{\cal B}}
\renewcommand{\P}{\mathbb P}
\newcommand{\E}{\mathbb E}
\newcommand{\dP}{P^{\downarrow}}
\newcommand{\dQ}{Q^{\downarrow}}
\newcommand{\couponm}{M}

\newcommand{\Nat}{\Bbb N}

\newcommand{\ed}{ \stackrel{d}{=}}

\def\endpf{\hfill $\Box$ \vskip0.5cm}
\def \proof{\noindent{\it Proof.\ }}

\newtheorem{theorem}{\large Theorem}
\newtheorem{proposition}[theorem]   {\large Proposition}

\newtheorem{corollary}[theorem]{\large  Corollary}
\newtheorem{lemma}[theorem]{\large  Lemma}

\begin{document}

\title{
Characterizations of exchangeable partitions and random discrete distributions by deletion properties \thanks{
This research is supported in parts by N.S.F.\
Award 0806118. }}
\author{
Alexander Gnedin%
\thanks{%
University of Utrecht; email A.V.Gnedin@uu.nl}
\and 
Chris Haulk%
\thanks{%
University of California at Berkeley; email haulk@stat.berkeley.edu} 
\and
Jim Pitman
\thanks{%
University of California at Berkeley; email pitman@stat.berkeley.edu} 
}
\date{\today}
\maketitle

\begin{abstract} 
\noindent
We prove a long-standing conjecture which characterises the Ewens-Pitman 
two-parameter family of exchangeable random partitions, plus a short list of limit
and exceptional cases, by the following property: for each $n = 2,3, \ldots$, if one of $n$ individuals is chosen uniformly at random,
 independently of the random partition $\pi_n$ of these individuals into various types, and all individuals of the same type as the
chosen individual are deleted, then for each $r > 0$, given that $r$ individuals remain, these individuals are partitioned 
according to $\pi_r'$ for some sequence of random partitions $(\pi_r')$
 which does not depend on $n$.
An analogous result characterizes the associated Poisson-Dirichlet family of random discrete distributions 
by an independence property 
related to random deletion of a frequency chosen by a size-biased pick. We also survey the regenerative properties of members of the two-parameter family, and settle a question regarding the explicit arrangement of intervals with lengths given by the terms of 
the Poisson-Dirichlet random sequence into the interval partition induced by the range of a homogeneous 
neutral-to-the right process.

\end{abstract}
\tableofcontents

\section{Introduction}

Kingman \cite{MR509954}
introduced the concept of a {\em partition structure}, that is a
family of probability distributions for random partitions $\pi_n$ of a positive integer $n$, with a sampling
consistency property as $n$ varies. 
Kingman's work was motivated by applications in population
genetics, where the partition of $n$ may be the allelic partition generated by randomly sampling a set of $n$ individuals from
a population of size $N \gg n$, considered in a large $N$ limit which implies
sampling consistency.  
Subsequent authors have established the importance of Kingman's theory of partition structures, 
and representations of these structures in terms of exchangeable random partitions and random discrete 
distributions \cite{jp.isbp}, in a number of other settings, which include the theory of species sampling 
\cite{jp96bl}, 
random trees and associated random processes of fragmentation and coalescence 
\cite{jp97cmc,hpw07,BertoinFragCoag,BertoinCoagFrag},
Bayesian statistics and machine learning 
\cite{teh06,09memoizer}.
Kingman 
\cite{MR0526801}
showed that the Ewens sampling formula from population genetics defines a particular partition 
structure $(\pi_n)$, which he characterized by the following property, together
with the regularity condition that $\P(\pi_n = \lambda) >0 $ for every
partition $\lambda$ of $n$:

\begin{quote}

for each $n = 2,3, \ldots$, if an individual is chosen uniformly at random 
independently
of a random partitioning of these individuals into various types according to $\pi_n$, and all individuals of the same type as the
chosen individual are deleted, then conditonally given that the number of remaining individuals is $r >0$, 
these individuals are partitioned according to a copy of $\pi_r$.  
\end{quote}
We establish here a conjecture of Pitman \cite{jp.epe}  that if this property is weakened by replacing $\pi_r$ by $\pi_r'$ 
for some sequence of 
random partitions $(\pi_r')$, and a suitable regularity condition is imposed, 
then 
$(\pi_n)$ belongs to the two-parameter family of partition structures
introduced in \cite{jp.epe}.
Theorem \ref{2pchar1} below provides a more careful statement.  
We also present a corollary of this result, to characterize 
the two-parameter family of Poisson-Dirichlet distributions by an independence property of a single size-biased pick,
thus improving upon \cite{jp.isbp}.

\par 
Kingman's characterization of the Ewens family of partition structures by 
deletion of a type has been extended in another direction  by allowing
other deletion algorithms but continuing to require that the distribution 
of the partition structure be preserved.
The resulting theory of {\em regenerative partition structures} \cite{rps04},
is connected to the theory  of regenerative sets, including Kingman's regenerative phenomenon \cite{KingmanReg}, on a multiplicative scale.
In the last section of the paper we review such deletion properties of the two-parameter family of partition structures, and offer a new proof of a result of
Pitman and Winkel \cite{PitmanWinkel} regarding the 
explicit arrangement of intervals with lengths given by the terms of the Poisson-Dirichlet random sequence into the interval partition induced by a
multiplicatively regenerative set.

\section{Partition Structures}

This section briefly reviews Kingman's theory of partition structures, which provides the general context of this article.
To establish some terminology and notation for use throughout the paper, recall
that a {\em composition } $\lambda$ of a positive integer $n$
is a sequence of positive integers
$\lambda = (\lambda_1, \ldots, \lambda_k)$, with $\sum_{i = 1}^{k} \lambda_i = n$.
Both $k = k_\lambda$ and $n = n_\lambda$ may be regarded as functions of $\lambda$. Each term $\lambda_i$ is called a 
{\em part} of $\lambda$. 
A {\em partition $\lambda$ of $n$} is a multiset of positive integers whose sum is $n$, commonly identified with the composition of $n$ obtained by putting its positive
integer parts in decreasing order, or with the infinite 
sequence of non-negative integers obtained by appending an infinite string of zeros  to this composition of $n$. So
$$\lambda = (\lambda_1, \lambda_2, \ldots) \mbox{~ with~ } \lambda_1 \ge  \lambda_2 \ge \cdots \ge 0  $$
represents a partition of $n = n_\lambda$ into $k = k_\lambda$ parts, where
$$
n_\lambda := \sum_i \lambda_i  \mbox{ and } k_\lambda:= \max \{ i : \lambda_i > 0 \}.
$$
Informally, 
a partition $\lambda$ describes an unordered collection of $n_\lambda$ balls of 
$k_\lambda$ different colors, with $\lambda_i$ balls of the $i$th most frequent color. 
A {\em random partition of $n$} is a random variable $\pi_n$ with values in the finite set of all partitions $\lambda$ of $n$.  
Kingman \cite{MR509954}
defined a {\em partition structure} to be a sequence of random partitions $(\pi_n)_{n \in \N}$ 
which is {\em sampling consistent} in the following sense:
\begin{quote}
if a ball is picked uniformly at random and deleted from $n$ balls randomly colored according to $\pi_n$, then
the random coloring of the remaining  $n-1$ balls is distributed according to $\pi_{n-1}$.
\end{quote}

As shown by Kingman 
\cite{MR671034},
the theory of partition structures and associated partition-valued processes is best developed in terms of random partitions of the
set of positive integers. Our treatment here follows \cite{jp.epe}.
If we regard a random partition $\pi_n$ of a positive integer $n$ as a random coloring of $n$ unordered balls, an associated 
random partition $\Pi_n$ of the set $[n]:= \{1, \ldots, n \}$ may be obtained by placement of the colored balls in a row.
We will assume for the rest of this introduction that this placement is made by a random permutation which given $\pi_n$ is uniformly distributed
over all $n!$ possible orderings of $n$ distinct balls.

Formally, a partition of $[n]$ is a collection of disjoint non-empty blocks
$\{B_1, \ldots, B_k\}$ with $\cup_{i=1}^k B_i = n$
for some $1 \le k \le n$, where each $B_i \subseteq [n]$ represents the set 
of places occupied by balls of some particular color.
We adopt the convention that the blocks $B_i$ are listed {\em in order of appearance}, meaning that $B_i$ is the set of places in the row occupied by balls 
of the $i$th color to appear.
So $1 \in B_1$, and if $k \ge 2$ the least element of $B_2$ is the 
least element of $[n]\setminus B_1$, if $k \ge 3$ the least element of $B_3$ is the least element of $[n]\setminus(B_1 \cup B_2)$, and so on.
This enumeration of blocks identifies each partition of $[n]$ with an {\it ordered partition} $(B_1,\ldots,B_k)$,
subject to these constraints. 
The sizes of parts $(|B_1|,\ldots, |B_k|)$ of this partition form a composition of $n$.
The notation $\Pi_n = (B_1, \ldots, B_k)$ is
used to signify that $\Pi_n = \{B_1, \ldots, B_k\}$ for some particular sequence of blocks 
$(B_1, \ldots, B_k)$ listed in order of appearance. 
If $\Pi_n$ is derived from $\pi_n$ by uniform random placement of balls in a row, then $\Pi_n$ is {\em exchangeable},
meaning that its distribution is invariant under every deterministic rearrangement of places by a permutation of $[n]$. 
Put another way,
for each partition $(B_1, \ldots, B_k)$ of $[n]$, with blocks in order of appearance,
\begin{equation}
\label{EPF}
\prob( \Pi_n = (B_1, \ldots, B_k) ) = p ( |B_1|, \ldots, |B_k|)
\end{equation}
for a function $p = p(\lambda)$ of compositions  $\lambda$ of $n$ which is a symmetric function of its $k$ arguments for each $1 \le k \le n$.
Then $p$ is called the {\em exchangeable partition probability function} (EPPF) associated with $\Pi_n$, or with $\pi_n$, the partition of $n$
defined by the unordered sizes of blocks of $\Pi_n$. 


As observed by Kingman 
\cite{MR671034}, 
$(\pi_n)$ is sampling consistent if and only if the sequence of partitions $(\Pi_n)$ can be
constructed to be consistent in the sense that for $m<n$ the restriction of $\Pi_n$ to $[m]$ is $\Pi_m$.  This amounts to a simple recursion
formula satisfied by $p$, recalled later as \re{add-rule}. 
The sequence $\Pi = (\Pi_n)$ can then be interpreted as a random partition of the set $\N$ of all positive integers, whose restriction to $[n]$ is $\Pi_n$ for every $n$.
Such $\Pi$ consists of a sequence of blocks $\B_1, \B_2, \ldots$, which may be identified as random disjoint subsets of $\N$, with
$\cup_{i = 1}^\infty \B_i = \N$, where the nonempty blocks 
are 
arranged by increase of their minimal elements, and 
if the number of nonempty blocks is some $K < \infty$, then by convention $\B_i = \varnothing$ for $ i > K$. 
Similarly, $\Pi_n$ consists of a sequence of blocks $\B_{ni}:=\B_{i}\cap [n]$, 
 where 
$\cup_i\B_{ni}=[n]$, and
the nonempty blocks are consistently arranged
by increase of their minimal elements,
for all $n$.

These considerations are summarized by the following proposition:
\begin{proposition}
\label{kinp}{\rm (Kingman 
\cite{MR671034}
)}
The most general partition structure, 
defined by a sampling consistent collection of distributions for 
partitions $\pi_n$ of integers $n$, is associated with a unique 
probability distribution of an exchangeable partition of positive integers $\Pi = (\Pi_n)$, as determined by an EPPF $p$ according to {\rm \re{EPF}}. 
\end{proposition}

We now recall a form of {\em Kingman's paintbox construction} of such an exchangeable random partition $\Pi$ of 
positive integers.
Regard the unit interval $[0,1)$ as a continuous spectrum of distinct colors, and suppose given 
a sequence of random variables $(\dP_1,\dP_2,\ldots)$ called
{\em ranked frequencies},
subject to  the constraints
\eq
\label{part1}
1\geq \dP_1\geq \dP_2\geq \ldots\geq0,~~~P_*:=1-\sum_{j=1}^\infty \dP_j\geq 0.
\en
The color spectrum is partitioned into a sequence of intervals $[l_i, r_i)$ of lengths $\dP_i$, and 
in case $P_* >0$ a further interval 
$[1-P_*, 1)$ of length $P_*$.  
Each point $u$ of $[0,1)$ is assigned the color $c(u) = l_i$ if $u \in [l_i, r_i)$ for some $i = 1,2, \ldots$,
and $c(u) = u$ if $u \in [1- P_*, 1)$. This coloring of points of $[0,1)$, called {\em Kingman's paintbox} 
associated with $(\dP_1,\dP_2,\ldots)$,
is sampled by an infinite sequence of independent uniform$[0,1]$ variables $U_i$, 
to assign a color $c(U_i)$ to the $i$th ball in a row of balls indexed by $i = 1,2, \ldots$.
The associated {\em color partition} of $\N$ is generated 
by the random equivalence relation $\sim$ defined by $i \sim j $ if and only if $c(U_i) = c(U_j)$,
meaning that either $U_i$ and $U_j$ fall in the same compartment  of the paintbox,  or that $i = j$ and $U_i$ falls in $[1-P_*,1)$.

\begin{theorem}\label{Kpaintbox}
{\em (Kingman's paintbox representation of exchangeable partitions 
\cite{MR509954})} 
Each exchangeable partition $\Pi$ of $\N$ generates a sequence of
ranked frequencies $(\dP_1,\dP_2,\ldots)$ such that the conditional distribution of $\Pi$  given these frequencies is that of the color partition of $\N$
derived from 
$(\dP_1,\dP_2,\ldots)$ by Kingman's paintbox construction.
The exchangeable partition probability function $p$ associated with $\Pi$ determines the distribution of $(\dP_1,\dP_2,\ldots)$, and vice versa.
\end{theorem}

The distributions of ranked frequencies $(\dP_1,\dP_2,\ldots, )$ associated with naturally arising partition structures
$(\pi_n)$ are quite difficult to deal with analytically. See for instance \cite{py95pd2}. Still, $(\dP_1,\dP_2,\ldots)$ can be constructed as the 
decreasing rearrangement of the frequencies $P_i$ of blocks $\B_i$ of $\Pi$
defined as the almost sure limits
\begin{equation}
\label{freqs}
P_i=\lim_{n\to\infty} n^{-1} |\B_{ni}|   
\end{equation}
where $i = 1,2, \ldots$ indexes the blocks in order of appearance, 
while
$$
P_* = 1 - \sum_{i = 1}^\infty P_i = 1 - \sum_{i = 1}^\infty \dP_i
$$
is the asymptotic frequency of the union of singleton blocks
$$\B_*:= \cup_{\{i : |\B_i| = 1\}} \B_i,$$
so that \re{freqs} holds also for $i = *$.
The frequencies are called {\em proper} if $P_* = 0$ a.s.; 
then almost surely every nonempty block $\B_i$ of $\Pi$ has a 
strictly positive frequency, hence $|\B_i| = \infty$, while every block $\B_i$ with $0 < |\B_i| < \infty$ is a singleton block.

The ranked frequencies $\dP_1, \dP_2, \ldots $ appear 
in the sequence $(P_j)$ in the order in which intervals of these lengths 
are discovered by a process of uniform random sampling, as in Kingman's paintbox construction. If $P_* >0$ then in addition to the strictly positive 
terms of $\dP_1, \dP_2, \ldots$ the sequence $(P_i)$ also contains infinitely many zeros which correspond to singletons in $\Pi$.  The conditional distribution of $(P_j$) given $(\dP_j)$ can also be described in terms of iteration of a single size-biased pick, defined as follows.  For a sequence of non-negative random variables $(X_i)$ with $\sum_i X_{i}\leq 1$ 
and a random index $J\in \{1,2,\ldots,\infty\}$, 
call $X_J$ a {\em size-biased pick} from $(X_i)$ if
$X_J$ has value $X_j$ if $J= j < \infty$ and $X_J = 0$ if $J = \infty$,
with
\begin{equation}
\label{sbpick}
\prob(J = j\, |\, (X_i, i\in \Nat)) = X_j
~~~( 0 < j < \infty )
\end{equation}
(see \cite{GnedinSBP} 
for this and another definition of size-biased pick in the case of improper frequencies).
The {\em sequence derived from $(X_i)$ by deletion
of $X_J$ and renormalization} refers to the sequence $(Y_i)$ obtained from $(X_i)$ by 
first deleting the $J$th term $X_J$, then closing up the gap if $J\neq\infty$, and finally normalizing each term by $1 - X_J$.
Here by convention, $(Y_i)$ = $(X_i)$ if $X_J= 0$ and $(Y_i)$ is the the zero sequence if $X_J=1$.
Then $P_1$ is a size-biased pick from $(\dP_j)$, $P_2$ is a size-biased pick from the sequence derived from $(\dP_j)$ by deletion of $P_1$ and renormalization, and so on.  For this reason, $(P_i)$ is said to be a \emph{size-biased permutation} of $(\dP_i)$.

\paragraph{The two-parameter family}
It was shown in \cite{jp.epe} that for each pair of real parameters $(\alpha,\theta)$ with
\eq\label{2Prange}
0\leq\alpha <1,~\theta>-\alpha
\en
the formula
\eq\label{EPPF2Par}
p_{\alpha,\theta} (\lambda) := {\prod_{i=1}^{k-1}(\theta+i\alpha)\over(\theta+1)_{n -1}}\prod_{j=1}^{k}(1-\alpha)_{\lambda_j - 1}
\en
where $k = k_\lambda$, $n = n_\lambda$, and
$$(x)_n := x (x+1) \ldots (x+n-1) = \frac{\Gamma(x +n)}{\Gamma(x)}$$
is a rising factorial, defines the EPPF of an exchangeable random partition of positive integers whose block frequencies $(P_i)$ in order
of appearance admit the {\em stick-breaking representation}
\begin{equation}\label{St-Br}
P_i=W_{i}\prod_{j=1}^{i-1} (1-W_j)
\end{equation}
for random variables $W_j$ such that
\begin{equation}
\label{wind2}
W_1, W_2, \ldots \mbox{ are mutually independent}
\end{equation}
with
\eq
\label{W-beta}
W_k\ed \beta_{1-\alpha,\theta+k\alpha}
\en
where  $\ed$ indicates equality in distribution, and $\beta_{a,b}$ for $a,b >0$ denotes a random variable
with the beta$(a,b)$ density 
\begin{equation}
\label{betadens}
\prob(\beta_{a,b}\in {\rm d}u ) = {\Gamma(a+b) \over \Gamma(a) \Gamma(b) } u ^{a-1}(1-u)^{b-1}{\rm d}u ~~~~~~(0 < u < 1)
\end{equation}
which is also characterized by the moments
\begin{equation}
\label{betamoms}
{\mathbb E} [\beta_{a,b}^i (1 - \beta_{a,b})^j ] =  \frac{ (a)_i  (b)_j }{ (a + b)_{i+j} } ~~~~~~~(i,j  = 0,1, 2, \ldots ).
\end{equation}

Formula \re{EPPF2Par} also defines an  EPPF for $(\alpha,\theta)$ in the range
\eq\label{alphaneg}
 \alpha<0,~\theta=-\couponm \alpha ~{\rm for~some~}\couponm  \in \Nat, 
\en
in which case the stick-breaking representation \re{St-Br} with factors as in \re{W-beta} makes sense for 
$1 \le  k \le \couponm $, with the last factor $W_\couponm  = 1$.  
The frequencies $(P_1, \ldots, P_\couponm)$ in this case are a size-biased random permutation of $(Q_1, \ldots, Q_\couponm)$ with
the symmetric Dirichlet distribution with $\couponm$ parameters equal to $\nu:= - \alpha >0$. It is well known that the $Q_i$ can be 
constructed as $Q_i = \gamma_{\nu}^{(i)}/\Sigma, 1 \le i \le \couponm$,
where $\Sigma = \sum_{i=1}^\couponm \gamma_\nu^{(i)}$ and the $\gamma_{\nu}^{(i)}$ are independent and identically distributed copies of a gamma variable $\gamma_\nu$  with density
\begin{equation}
\label{gammadens}
\P( \gamma_\nu \in {\rm d}x ) = \Gamma(\nu)^{-1} x ^{\nu-1} e^{-x} {\rm d}x  ~~~~~~( x > 0 ).
\end{equation}
As shown by Kingman \cite{MR0526801},
the $(0,\theta)$ EPPF \re{EPPF2Par} for $\alpha = 0, \theta >0$ arises in the limit of random sampling from such symmetric Dirichlet frequencies as $\nu = -\alpha \downarrow 0$ and $\couponm \uparrow \infty$ with
$\nu \couponm = \theta$ held fixed. In this case, the distribution of the partition $\pi_n$ is that determined by the Ewens sampling formula with parameter $\theta$, the residual fractions
$W_i$ in the stick-breaking representation are identically distributed like $\beta_{1,\theta}$, and the ranked frequencies $\dP_i$ can be obtained by normalization of the jumps of a
gamma process with stationary independent increments $(\gamma_\nu,  0 \le \nu \le \theta)$. Perman, Pitman and Yor \cite{ppy92} gave extensions of this description to the case 
$0 < \alpha < 1$ when the distribution of ranked frequencies can be derived from the jumps of a stable subordinator of index $\alpha$.  
See also \cite{py95pd2,pitman02pk,csp} for further discussion and applications to the description of ranked lengths of excursion intervals of Brownian motion and Bessel processes.

In the limit case when $\nu = - \alpha \to \infty$ and $\theta = \couponm \nu \to \infty$, for a fixed positive integer $\couponm$, the EPPF \re{EPPF2Par}
converges to 
\eq\label{Coupon}
p_{\couponm}
(\lambda):= \frac{ \couponm (\couponm -1) \cdots (\couponm -k+1) }{\couponm ^n}\,,
\en
corresponding to sampling from $\couponm$ equal frequencies
$$
P_1 = P_2 = \cdots = P_\couponm = 1/\couponm 
$$
as in the classical coupon collector's problem with some fixed number $\couponm$ of equally frequent types of coupon.
We refer to the collection of partition structures defined by \re{EPPF2Par} for the parameter
ranges (\ref{2Prange}) and (\ref{alphaneg}), 
as well as the limit cases \re{Coupon}, as the 
{\em extended two-parameter family}.

\par The partition $\mathbf{0}$ of $\mathbb{N}$ into singletons and the partition $\mathbf{1}$ of $\mathbb{N}$ into a single block 
both belong to the closure of the two-parameter family.
As noticed by Kerov \cite{Kerov}, a mixture of these two trivial partitions with mixing proportions $t$ and  $1-t$ 
also belongs to the closure, as is seen 
from (\ref{EPPF2Par}) by letting $\alpha \to 1$ and $\theta \to -1$ in such a way that 
$(1- \alpha)/(\theta+1)\to t\,$ and $(\theta+\alpha)/(\theta +1)\to 1-t$.


\paragraph{Characterizations by deletion properties}
The main focus of this paper is the following result, which was first conjectured by Pitman \cite{jp.epe}.
For convenience in presenting this result, we impose the following mild {\em regularity condition}
on the EPPF $p$ associated with a partition structure $(\pi_n)$:
\eq
\label{p221}
p(2,2,1) > 0  \mbox{~ and~ } \lim_{n \te \infty}p(n) = 0.
\en
Equivalently, in terms of the frequencies $P_i$ in order of appearance,
\eq
\label{p221x}
\prob(0 < P_1 < P_1 + P_2 < 1 ) > 0  \mbox{ and } \prob(P_1 = 1) = 0,
\en
or again, in terms of the ranked frequencies $\dP_i$,
\eq
\label{p221r}
\prob(0 < \dP_2 , ~\dP_1 + \dP_2 < 1 ) > 0 \mbox{ and } \prob(\dP_1 = 1) = 0.
\en
Note that this regularity condition does not rule out the case of improper frequencies.
See Section \ref{secfails} for discussion of how the following results can be modified to accomodate partition 
structures not satisfying the regularity condition.

\begin{theorem}\label{2pchar1}
Among all partition structures $(\pi_n)$
with
EPPF $p$ subject to {\rm \re{p221}},
the extended two-parameter family is characterized by the following property:
\begin{quote}
if one of $n$ balls is chosen uniformly at random, independently of a random coloring of these balls according to $\pi_n$, then 
given the number of other balls of the same color as the chosen ball is $m-1$, for some $1 \le m < n$, the coloring of 
the remaining $n-m$ balls is distributed according to $\pi_{n-m}'$ for some sequence of partitions
$(\pi_1',\pi'_2,\ldots)$ which does not depend on $n$.
\end{quote}
Moreover, if $(\pi_n)$ has the $(\alpha,\theta)$ EPPF {\rm\re{EPPF2Par}},  
then $(\pi_n')$ has the $(\alpha,\theta +\alpha )$ EPPF {\rm\re{EPPF2Par}},
whereas if $(\pi_n)$ has the EPPF {\rm \re{Coupon}} for some $\couponm$, 
then the EPPF of $(\pi_n')$ has the same form except with $\couponm$ decremented by $1$.
\end{theorem}
\noindent
Note that it is not assumed as part of the property that $(\pi_n')$ is a 
partition structure. Rather, this is implied by the conclusion.  
Our formulation of Theorem \ref{2pchar1} was inspired by Kingman 
\cite{MR0526801}
who assumed also that $\pi_n' \stackrel{d}{=} \pi_n$ for all $n$. The conclusion then holds with $\alpha = 0$, in which case the distribution of
$\pi_n$ is that determined by the Ewens sampling formula from population genetics.  

In Section \ref{sect:ExchPart} we offer a proof of Theorem \ref{2pchar1} by purely combinatorial methods.
Some preliminary results which we develop in Section \ref{sect:PEP} allow
Theorem \ref{2pchar1} to be reformulated in terms of frequencies as in the following Corollary:

\begin{corollary}\label{2pcharcor}
Let the asymptotic frequencies $(P_i)$ of an exchangeable random partition of positive integers $\Pi$ be represented in the
stick-breaking form {\rm \re{St-Br}} for some sequence of random variables $W_1<1, W_2, \ldots$.
The condition
\begin{equation}
\label{wind1}
W_1 \mbox{ is independent of } (W_2, W_3, \ldots )
\end{equation}
obtains if and only if 
\begin{equation}
\label{windm}
\mbox{ the $W_i$ are mutually independent. }
\end{equation}
If in addition to $(\ref{wind1})$ the regularity condition {\rm \re{p221}} holds, then
$\Pi$ is governed  by the extended two-parameter family, either with $W_i \ed \beta_{1-\alpha,\theta+i\alpha}$, or with $W_i = 1/(\couponm - i +1)$ for $1 \le i \le \couponm$, as in the limit case {\rm \re{Coupon}}, for some 
$\couponm  = 3,4, \ldots$.  
\end{corollary}
\noindent
The characterization of the two-parameter family using \re{windm} rather than the weaker condition \re{wind1} 
was provided by Pitman \cite{jp.isbp}.  
As we show in Section \ref{sect:ExchPart}, it is possible to derive \re{windm} directly from \re{wind1}, without passing via Theorem \ref{2pchar1}.

The law of frequencies $(P_i)$ defined by the stick-breaking scheme
\re{St-Br} for independent factors $W_i$ with
$W_i \ed \beta_{1-\alpha,\theta+i\alpha}$ is
known as the the two-parameter Griffiths-Engen-McCloskey distribution, 
denoted ${\rm GEM}(\alpha,\theta)$.
The property of the independence of residual proportions $W_i$, also known as 
{\it complete neutrality}, has also been studied extensively in connection 
with finite-dimensional Dirichlet distributions \cite{BobWes}.

The above results can also be expressed in terms of ranked frequencies. Recall that the distribution of ranked frequencies $(\dP_k)$ of an $(\alpha,\theta)$-partition is known as the
{\em two-parameter Poisson-Dirichlet distribution ${\rm PD}(\alpha,\theta)$}.  
According the the previous discussion, a random sequence $(\dP_k)$ with 
${\rm PD}(\alpha,\theta)$ distribution is obtained by ranking 
a sequence $(P_i)$ with ${\rm GEM}(\alpha,\theta)$ distribution.
The ${\rm PD}(\alpha,\theta)$ distribution was systematically studied in \cite{py95pd2},
and has found numerous further applications to random trees and associated processes of fragmentation and coagulation 
\cite{jp97cmc,hpw07,BertoinCoagFrag}.

\begin{corollary}\label{corextra}
Let $(\dP_k)$ be a decreasing sequence of ranked frequencies subject to the regularity condition {\rm(\ref{part1})} and {\rm(\ref{p221r})}.
For $\dP_J$,  a size-biased pick from $(\dP_k)$, let $(\dQ_k)$ be derived from
$(\dP_k)$ by deletion of $\dP_J$ and renormalization.  The random variable $\dP_J$ is independent of the sequence $(\dQ_k)$ if and only if either the distribution of
$(\dP_k)$ is ${\rm PD}(\alpha,\theta)$ for some $(\alpha,\theta)$, or $\dP_k = 1/\couponm $ for all $1 \le k \le \couponm $, for some $\couponm  
\ge 3$. In the former case, the distribution of 
$(\dQ_k)$ is ${\rm PD}(\alpha,\theta + \alpha)$, 
whereas in the latter case, the deletion and renormalization simply decrements $\couponm$ by one.
\end{corollary}
The `if' part of this Corollary is Proposition 34 of Pitman-Yor \cite{py95pd2}, while the `only if' part follows
easily from Corollary \ref{2pcharcor}, using Kingman's paintbox representation.

\section{Partially Exchangeable Partitions}\label{sect:PEP}

We start by recalling from \cite{jp.epe} some basic properties of {\em partially exchangable partitions of positive
integers}, which are consistent sequences $\Pi = (\Pi_n)$, where $\Pi_n$ is a partition of $[n]$ whose probability
distribution is of the form \re{EPF} for some function $p = p(\lambda)$ of compositions $\lambda$ of positive
integers. The consistency of $\Pi_n$ as $n$ varies amounts to the {\em addition rule}
\eq\label{add-rule}
p(\lambda)=\sum_{j =1}^{k + 1 } p(\lambda^{(j)}),
\en
where $k = k_\lambda$ is the number of parts of $\lambda$, and $\lambda^{(j)}$  is the composition of $n_\lambda +1$ derived from $\lambda$ by incrementing $\lambda_j$ to $\lambda_j+1$,
and leaving all other components of $\lambda$ fixed. In particular, 
for $j = k_\lambda +1$ this means appending a $1$ to $\lambda$.
There is also the normalization condition $p(1)=1$. To illustrate \re{add-rule} for $\lambda = (3,1,2)$:
$$
p(3,1,2) = p(4,1,2) + p(3,2,2) + p(3,1,3) + p(3,1,2,1).
$$
The following proposition recalls the analog of Kingman's representation
for partially exchangeable partitions:

\begin{proposition}\label{pep} {\rm (Corollary 7 from \cite{jp.epe})}
Every partially exchangeable partition of positive integers
$\Pi$ is such that for each $k \ge 1$, the $k$th block ${\cal B}_k$ has an almost sure limit frequency $P_k$.
The partition probability function $p$ can then be presented as
\begin{equation}\label{Eqn:EPPFMomentFormula}
p(\lambda)
=
{\mathbb E}\left[ \prod_{i=1}^k 
P_i^{\lambda_i-1}\prod_{j=1}^{k-1} R_j \right] ,
\end{equation}
where $k =k_\lambda$ and $R_j := (1 - P_1 - \cdots -P_j)$.
Alternatively, in terms of 
the residual fractions $W_k$ in the stick-breaking representation {\rm \re{St-Br}}:
\begin{equation}\label{Eqn:EPPFMomentFormula2}
p(\lambda) = {\mathbb E}\left[\prod_{i=1}^k W_i^{\lambda_{i}-1}\overline{W}_i^{\Lambda_{i+1}}\right], 
\end{equation}
where $\overline{W}_i:= 1- W_i$, $\Lambda_j:=\sum_{i\geq j}\lambda_i$. 
This formula sets up a correspondence between the probability distribution of $\Pi$, 
encoded by the partition probability function $p$,
and an arbitrary joint distribution of a sequence of random variables 
$(W_1, W_2, \ldots)$ with $0 \le W_i \le 1$ for all $i$.
\end{proposition}
In terms of randomly coloring a row of $n_\lambda$ balls, the product whose expectation appears in \re{Eqn:EPPFMomentFormula2} is the conditional probability given $W_1, W_2, \ldots$ of the event that
the first $\lambda_1$ balls are colored one color, the next $\lambda_2$ balls  another color, and so on. So \re{Eqn:EPPFMomentFormula2} reflects the fact that conditionally given 
$W_1, W_2, \ldots$ the process of random coloring of integers occurs according to the following {\rm residual allocation scheme}
\cite[Construction 16]{jp.epe}: 
\begin{quote}
Ball $1$ is painted a first color, 
and so is each subsequent ball according to a sequence of independent trials with probability $W_1$ of painting with color 1.  
The set of balls so painted defines the first block $\B_1$ of $\Pi$.
Conditionally given $\B_1$, the first unpainted ball is painted a second color, 
and so is each subsequent unpainted ball according to a sequence of independent trials with probability $W_2$ of painting with color 2.  
The balls colored 2 define $\B_2$, and so on.
Given an arbitrary sequence of random variables $(W_k)$ with $0 \le W_k \le 1$, this coloring scheme shows how to construct a partially exchangeable partition of $\N$ whose asymptotic block
frequencies are given by the stick-breaking scheme \re{St-Br}.
\end{quote}
Note that the residual allocation scheme terminates at the first $k$, if any, such that $W_k = 1$, by painting all remaining balls color $k$. 
The values of $W_i$ for $i$ larger than such a $k$ have no effect on the construction of $\Pi$, so cannot be recovered from its almost sure limit frequencies.
To ensure that a unique joint distribution of $(W_1, W_2, \ldots)$ is associated with each $p$,
the convention may be adopted that the sequence $(W_i)$ terminates at the first $k$ if any such that $W_k = 1$.
This convention will be adopted in the following discussion.

For $W_i$ which are independent, formula (\ref{Eqn:EPPFMomentFormula2}) factorizes as
\begin{equation}\label{Eqn:EPPFMomentFormula22}
p(\lambda) = \prod_{i=1}^k 
\E( W_i^{\lambda_{i}-1}\overline{W}_i^{\Lambda_{i+1}} ) .
\end{equation}
In particular, for independent $W_i$ with the beta distributions (\ref{W-beta}),  
this formula is readily evaluated using \re{betamoms}
to obtain \re{EPPF2Par}. Inspection of \re{EPPF2Par} shows that this function of compositions $\lambda$ is a symmetric function of its parts. Hence the
associated random partition $\Pi$ is exchangeable. 

There is an alternate sequential construction of the two-parameter family of partitions which has become known as the ``Chinese Restaurant Process''
(see \cite{csp}, {Chapter 3}). 
Instead of coloring rows  of balls,
imagine customers entering a restaurant with an unlimited number of tables. Initially customer $1$ sits at table $1$. At stage $n$, if there are $k$ occupied tables, the $i$th of them  occupied by $\lambda_i$
customers for $1 \le i \le k$,  customer $n+ 1$ sits at one of the
previously occupied tables with probability $(\lambda_i - \alpha)/(n + \theta)$, and  occupies a new table $k+1$ with probability 
$(\theta + k\alpha)/(n + \theta)$.
It is then readily checked that for each partition of $[n]$ into blocks $B_i$ with $|B_i| = \lambda_i$, after $n$ customers labeled by $[n]$ have entered the restaurant, the probability 
that those customers
labeled by $B_i$ sat at table $i$ for each $1 \le i \le k_\lambda$ is given by the product formula \re{EPPF2Par}. 
Moreover, the stick-breaking description of the limit frequencies $P_i$ is readily derived 
from the P{\'o}lya urn-scheme description of exchangeable trials which given a beta$(a,b)$-distributed variable $S$,
 are independent with success probability $S$.

Continuing the consideration of a partially exchangeable partition $\Pi$ of positive integers, we record the following Lemma.

\begin{lemma}\label{p.prime}
Let $\Pi$ be a partially exchangeable random partition of $\N$ with partition probability function $p$, and with blocks $\B_1, \B_2, \ldots$ 
and residual frequencies $W_1, W_2, \ldots$
such that $W_1 < 1$ almost surely.
Let $\Pi '$ denote the partition of $\N$ derived from $\Pi$ by deletion of the block $\B_1$ containing $1$ and re-labeling of $\N - \B_1$ by the increasing bijection with $\mathbb N$.
Then the following hold:
\begin{itemize}
\item[\rm(i)] The partition $\Pi'$ is partially exchangeable, with partition probability function
\eq
\label{piprime}
p'(\lambda_2,\ldots,\lambda_k)=\sum_{\lambda_1=1}^\infty {\lambda_1+\ldots+\lambda_k-2\choose \lambda_1-1}p(\lambda_1,\lambda_2,\ldots,\lambda_k)
\en
and residual frequencies $W_2, W_3, \ldots$.
\item[\rm(ii)] If $\Pi$ is exchangeable, then so is $\Pi'$.
\item[\rm(iii)] For $1 \le m \le n$ 
\eq
\label{qnm}
q(n:m) := \prob( \B_1\cap [n]=[m]) = {\mathbb E}(W_1^{m-1}\overline{W}_1^{n-m}),
\en
and there is the addition rule
\eq\label{add-rule-q}
q(n:m)=q(n+1:m+1)+q(n+1:m).
\en
\item[\rm(iv)]  Let $T_n := \inf \{m : |[n+m]\setminus\B_1| = n\}$ which is the number of balls of the first
color preceding the $n$th ball not of the first color. Then 
\eq
\label{negbin}
\prob (T_n = m) = {m+n-2\choose m-1} q(n+m:m),
\en
and consequently 
\begin{equation}\label{sumone}
\sum_{m=1}^\infty  {m+n-2\choose m-1} q(n+m:m) =1.
\end{equation}
\end{itemize}
\end{lemma}
\proof
Formula \re{qnm} is read from the general construction of $\B_1$ given $W_1$ by assigning each $i \ge 2$
to $\B_1$ independently with the same probability $W_1$.
The formulas \re{piprime} and \re{negbin} are then seen to be marginalizations of the following
expression for the joint distribution of $T_n$ and $\Pi_n'$, the restriction of $\Pi'$ to $[n]$:
\eq
\label{tnprime}
\prob( T_n = m, \Pi_n' = (C_1, \ldots, C_{k-1} ) )
= {m+n-2\choose m-1} q(n+m:m) p(m, |C_1|, \ldots, |C_{k-1}|)
\en
for every partition $(C_1, \ldots, C_{k-1})$ of $[n]$. To check \re{tnprime}, observe that the
event in question occurs if and only if $\Pi_{n+m} = (B_1, \ldots, B_{k})$ for some blocks $B_i$ with
$|B_1| = m$ and $|B_i| = |C_{i-1}|$ for $2 \le i \le k$. Once $B_1$ is chosen, each $B_i$ for 
$2 \le i \le k$ is the image of $C_{i-1}$ via the increasing bijection from $[n]$ to $[n+m] \setminus B_1$.
For prescribed $C_{i-1}, 2 \le i \le k$, the choice of $B_1 \subset [n+m]$ is arbitrary subject to the constraint that
$1 \in B_1$ and $n+m \notin B_1$. The number of choices is the binomial coefficient in \re{tnprime}, so the
conclusion is evident.
\endpf

The connection between Theorem \ref{2pchar1}  and Corollary \ref{2pcharcor}  is established by the following Lemma:

\begin{lemma}
\label{lemma4equiv}
Let $\Pi$ be a partially exchangeable partition of $\N$ with residual frequencies $W_i$ such that
$\prob (W_1 < 1 ) = 1$, with the convention that the sequence terminates at the first $k$ (if any) such that $W_k = 1$,
so the joint distribution of $(W_i)$ is determined uniquely by the partition probability function $p$ of $\Pi$,
and vice versa, according to formula 
{\rm \re{Eqn:EPPFMomentFormula2}}.
For $\B_1$  the first block of $\Pi$ with frequency $W_1$,
let $\Pi'$  be derived from $ \Pi$ by deleting block $\B_1$ 
and relabeling the remaining elements 
as in {\rm Lemma \ref{p.prime}}.
The following four conditions on $\Pi$ are equivalent:
\begin{itemize}
\item[\rm(i)] $W_1$ is independent of $(W_2, W_3, \ldots)$.
\item[\rm(ii)]  The partition probability function $p$ of $\Pi$ admits 
a factorization of the following form, for all compositions $\lambda$ of positive integers with $k \ge 2$ parts:
\eq
\label{factor}
p(\lambda) = q( n_\lambda: \lambda_1) p'(\lambda_2, \ldots, \lambda_k)
\en
for some non-negative functions $q(n:m)$ and $p'(\lambda_2, \ldots, \lambda_k)$.
\item[\rm(iii)]  
For each $1 \le m < n$, the conditional distribution of $\Pi_{n-m}'$ given $|\B_1 \cap [n]| = m$ depends only on $n-m$.
\item[\rm(iv)]  
The random set $\B_1$ is independent of the random partition $\Pi'$ of $\N$.
\end{itemize}
Finally, if these conditions hold, then {\rm (ii)} holds in particular for $q(n:m)$ as in {\rm \re{qnm}} 
and $p'(\lambda_2, \ldots, \lambda_k)$ 
the partition probability function of $\Pi'$.
\end{lemma}
\proof
That (i) implies (ii) is immediate by combination of the moment formula
\re{Eqn:EPPFMomentFormula2}, \re{piprime} and
\re{qnm}.
Conversely, if (ii) holds for some $q(n:m)$ and $p'(\lambda_2, \ldots, \lambda_k)$, Lemma \ref{p.prime} 
implies easily that (ii) holds for $q$ and $p'$ as in that
Lemma. So (ii) gives a formula of the form
\eq
\label{fg}
{\mathbb E} [ f(W_1) g( W_2, W_3, \ldots ) ] = {\mathbb E} [ f(W_1) ] {\mathbb E} [ g(W_2, W_3, \ldots )],
\en
where $g$ ranges over a collection of bounded measurable functions whose expectations determine the law of
$W_2, W_3, \ldots$, and for the $g$ associated with $\lambda_2, \ldots, \lambda_k$, the function $f(w)$
ranges over the polynomials $w^{m-1} (1-w)^n$ where $m = \lambda_1 \in \N$
and $n = n_\lambda -\lambda_1 = \sum_{j=2}^k \lambda_j$. But linear combinations of these polynomials can be used to
uniformly approximate any bounded continuous function of $w$ on $[0,1]$ which vanishes in a neighbourhood of $1$.
It follows that \re{fg} holds for all such $f$, for each $g$, hence the full independence condition (i).
Lastly, the equivalence of (ii), (iii) and (iv) is easily verified.
\endpf

\section{Exchangeable Partitions}\label{sect:ExchPart}

For a block $B$ of a random partition $\Pi_n$ of $[n]$ with $|B| = m$, let $\Pi_n \setminus B$ denote the partition of $[n-m]$ obtained by first deleting the 
block $B$ of $\Pi_n$, then mapping the restriction of $\Pi_n$ to $[n] \setminus B$ to a partition of $[n-m]$ via the increasing bijection
between $[n] \setminus B$ and $[n-m]$. In terms of a coloring of $n$ balls in a row, this means deleting all $m$ balls
of some color, then closing up the gaps between remaining balls, to obtain a coloring of $n-m$ balls in a row. 
Theorem \ref{2pchar1} can be formulated a little more sharply as follows:

\begin{theorem}\label{Kchar1}
Among all exchangeable partitions $(\Pi_n)$ of positive integers with EPPF $p$ subject to {\rm\re{p221}}, 
the extended two-parameter family is characterized by the following property:
\begin{quote}
if $\B_{n1}$ denotes the random block of $\Pi_n$ 
containing $1$, then for each $1 \le m < n$, 
conditionally given $\B_{n1}$ with $|\B_{n1}| = m$, 
the partition $\Pi_n \setminus \B_{n1}$ has the same distribution as $\Pi_{n-m}'$
for some sequence of partitions $\Pi'_1,\Pi'_2,\ldots$ which does not depend on $n$.
\end{quote}
Moreover, if $(\Pi_n)$ is an $(\alpha,\theta)$ partition, then we can take for $(\Pi_n')$ 
the exchangeable $(\alpha,\theta +\alpha )$ partition of $\Nat$. 
\end{theorem}

For an arbitrary partition $\Pi_n$ of $[n]$ with blocks listed in the order of appearance,
define  $J_n$ as the index of the block containing an element chosen from $[n]$ uniformly 
at random, independently of $\Pi_n$. 
We call the block  ${\cal B}_{nJ_n}$ a {\it size-biased pick} from the sequence of blocks.
 Note that this definition agrees with (\ref{sbpick}) in the sense that the number
$|\B_{nJ_n}|/n$ is a size-biased pick from the numerical sequence $(|{\cal B}_{nj}|/n, ~j=1,2,\dots)$, 
because 
given a sequence of blocks of partition $\Pi_n$ 
the value $J_n=j$ is taken with probability $|{\mathcal B}_{nj}|/n$.
Assuming $\Pi_n$ exchangeable, the size of the block $|{\mathcal B}_{n1}|$ 
has the same distribution as $|\B_{nJ_n}|$ conditionally given the ranked sequence of block-sizes,
and the reduced partitions $\Pi_n\setminus\B_{n1}$ and $\Pi_n\setminus\B_{nJ_n}$ 
also have the same distributions.  
The equivalence of Theorem \ref{2pchar1} and Theorem \ref{Kchar1} is evident from these considerations.


We turn to the proof of Theorem \ref{Kchar1}.
The condition considered in Theorem \ref{Kchar1} is just that considered in Lemma \ref{lemma4equiv}(iii),
so we can work with the equivalent factorization condition \re{factor}.
We now invoke the symmetry of the EPPF for an exchangeable $\Pi$.
Suppose that an EPPF $p$ admits the factorization \re{factor}, and re-write 
the identity \re{factor} in the form
$$p(m,\lambda)={q(|\lambda|+m:m)\over q(|\lambda|+1:1)} p(1,\lambda).$$
For this expression we must have non-zero denominator, but this is assured by $\P(0 < W_1 < 1) >0$, which
is implied by the regularity condition (\ref{p221}).
Instead of part  $m$ in $p(m,\lambda)$, we have now  $1$ in $p(1,\lambda)$.
But $p$ is symmetric, hence
we can  iterate, eventually reducing each part to $1$.
\par
Let $\lambda=(\lambda_1,\ldots,\lambda_k)$ be a generic composition, and denote $\Lambda_j=\lambda_j+\cdots+\lambda_k$ the tail 
sums, thus $\Lambda_1=|\lambda|$. Iteration yields
\eq\label{iterate}
p(\lambda)={q(\Lambda_1:\lambda_1)\over q(1+\Lambda_2:1)}{q(1+\Lambda_2:\lambda_2)\over q(2+\Lambda_3:1)}\cdots
{q(k-2+\Lambda_{k-1}:\lambda_{k-1})\over q(k-1+\Lambda_k:1)}
{q(k-1+\Lambda_{k}:\lambda_{k})\over q(k:1)}\,p(1^k),
\en
where $p(1^k)$ is the probability of the singleton partition of $[k]$.
This leads to the following lemma, which is a simplification of \cite[Lemma 12]{jp.isbp}:

\begin{lemma} \label{beta}
Suppose that an EPPF $p$ satisfies the factorization condition {\rm (\ref{factor})} and 
the regularity condition {\rm (\ref{p221})}.
Then  
\begin{itemize}
\item[{\rm (i)}]
either 
$$q(n:m)={(a)_{m-1}(b)_{n-m}\over (a+b)_{n-1}}$$ 
for some $a,b>0$, corresponding to $W_1$ with beta$(a,b)$ distribution,
\item[{\rm (ii)}] or 
$$q(n:m)= c^{m-1}(1-c)^{n-m}$$ 
for some $0<c<1$, corresponding to $W_1 = c$, in which case necessarily $c = 1/M$ for some $M \ge 3$.
\end{itemize}
\end{lemma}
\proof
By symmetry and the assumption that $p(2,2,1) >0$, it is easily seen from Kingman's paintbox representation that
for each $m = 1,2, \ldots$ there is some composition $\mu$ of $m$ such that
$$
p(3,2,\mu)=p(2,3,\mu)>0,
$$
where for instance $(3,2,\mu)$ means the composition of $5 + m$ obtained by concatenation of $(3,2)$ and $\mu$.
Indeed, it is clear that one can take either $\mu = 1^m$ or $\mu$ to be a single part of size $m$, according to
whether the probability of at least three non-zero frequencies is zero or greater than zero.
Applying (\ref{iterate}) for suitable $k \ge 3$ with $p(1^k) >0$, and cancelling some common factors of
the form $q(n',m')$, which are all strictly positive because $p(2,2,1) >0$ implies $\P(0 < W_1 < 1) >0$, we
see that for every $m = 1,2, \ldots$
\begin{equation}\label{rec-start}
{q(m+5:3)q(m+3:2)\over q(m+3:1)}={q(m+5:2)q(m+4:3)\over q(m+4:1)}.
\end{equation}
We have by the addition rule (\ref{add-rule-q})
$$q(m+1:2)=q(m:1)-q(m+1:1),~~~q(m+2:3)=q(m:1)-2q(m+1:1)+q(m+2:1),$$
and introducing variables
$x_m=q(m:1)$, $n=m+2$
$$ {(x_{n+1}-2x_{n+2}+x_{n+3})(x_n-x_{n+1})\over x_{n+1}}={(x_{n+2}-x_{n+3})(x_n-2x_{n+1}+x_{n+2})\over x_{n+2}}.$$
The recursion is homogeneous, to pass to inhomogeneous variables 
divide both sides of the equality by $x_n$, then
set $y_n:=x_{n+1}/x_n$ and rewrite as
$$(1-2y_{n+1}+y_{n+2}y_{n+1})(1-y_n)=(1-y_{n+2})(1-2y_n +y_ny_{n+1}),$$
which simplifies as
$$-2y_{n+1}+y_{n+1}y_{n+2}+y_ny_{n+1}=-y_n-y_{n+2}+2y_n y_{n+2}.$$
Finally, 
use substitution 
$$y_n=1-{1\over z_n}$$
to arrive at 
$${z_n-2 z_{n+1}+z_{n+2}\over z_nz_{n+1}z_{n+2}}=0.$$
From this, $z_n$ is a linear function of $n$, which must be nondecreasing to agree with $0<y_n<1$. 
\par If $z_n$ is not constant, then going back to $x_n$'s we obtain 
$$q(n:1)=  c_0 {(b)_{n-1}\over (a+b)_{n-1}}, ~~~n\geq 3,$$
for some $a,b,c_0$,  where the factor $c_0$ appears since the relation (\ref{rec-start}) is homogeneous.
It is seen from the moments representation
$$q(n:1)=\int_{[0,1]} (1-x)^{n-1}{\mathbb P}({P}_1\in {\rm d}x),~~~~n\geq 3,$$ 
that when $a,b$ are fixed, the factor $c_0$ is determined from the normalization by choosing a value of ${\mathbb P}({P}_1=1)$.
The condition $p(n)\to 0$ means that  ${\mathbb P}({P}_1=1)=0$, in which case  
$c_0=1$ and the distribution of ${P}_1$ is beta$(a,b)$ with some positive $a,b$.

\par If $(z_n, ~n\geq 3)$ is a  constant sequence, then
$q(n:1)$ is a geometric progression, and a similar argument shows that the case (ii) prevails. 
That $c = 1/\couponm$ for some $\couponm \ge 3$ is quite obvious: the only way that a size-biased choice of a frequency can be constant
is if there are $\couponm$ equal frequencies for some $\couponm \ge 1$. The regularity assumption 
\re{p221} rules out the cases $\couponm =1,2$.

\endpf
\noindent

\paragraph{Proof of Theorem \ref{Kchar1}}
In the case (i) of Lemma \ref{beta}, substituting in (\ref{iterate}) yields
$$\frac{p(\lambda)}{p(1^k)} ={(a)_{\lambda_1-1}(b)_{\Lambda_2}\over(a+b)_{\Lambda_1-1}}{(a+b)_{\Lambda_2}\over(b)_{\Lambda_2}}
{(a)_{\lambda_2-1}(b)_{\Lambda_3+1}\over(a+b)_{\Lambda_2}}{(a+b)_{\Lambda_3+1}\over(b)_{\Lambda_3+1}}\cdots
{(a)_{\lambda_k-1}(b)_{k-1}\over (a+b)_{\Lambda_k+k-2}}{(a+b)_{k-1}\over (b)_{k-1}},
$$
provided $p(1^k) >0$.  After cancellation this becomes
$$\frac{p(\lambda)}{p(1^k)} ={(a+b)_{k-1} \over (a+b)_{n-1}}\prod_{j=1}^k (a)_{\lambda_j-1},$$
where $n=\Lambda_1=\lambda_1+\ldots+\lambda_k=|\lambda|$. 
Specializing,
$${p(2,1^{k-1})\over p(1^k)}={a\over a+b+k-1}$$
and using the addition rule \re{add-rule}
$$p(1^k)=k p(2,1^{k-1})+p(1^{k+1}),$$
we obtain the recursion
$${p(1^{k+1})\over p(1^k)}={a+b+k(1-a)-1\over a+b+k-1},~~~p(1)=1.$$
Now (\ref{EPPF2Par}) follows readily by re-parametrisation $\theta=a+b-1,\,\alpha=1-a$.
\par The case (ii) of Lemma \ref{beta} is even simpler, as it is immediate that $W_1 = 1/\couponm$ implies that the
partition is generated as if by coupon collecting with $M$ equally frequent coupons.
\endpf

\paragraph{Proof of Corollary \ref{2pcharcor}}

As observed earlier, Corollary \ref{2pcharcor} characterizing the extended two-parameter family by the condition that
\eq
\label{w123}
W_1 \mbox{ and } (W_2, W_3, \ldots) \mbox{ are independent}
\en
can be read from Theorem \ref{Kchar1} and Lemma \ref{lemma4equiv}.
We find it interesting nonetheless to provide another proof of Corollary \ref{2pcharcor} based on analysis of the limit frequencies rather than the EPPF.
This was in fact the first argument we found, without which we might not have persisted with the algebraic approach of the previous section.

Suppose then that $W_1,W_2, W_3, \ldots$ is the sequence of residual fractions associated with an EPPF $p$, and that \re{w123} holds.
The symmetry condition $p(r+1,s+1) = p(s+1,r+1)$
and the moment formula \re{Eqn:EPPFMomentFormula2} give
\eq
\label{momeq}
{\mathbb E}(W_1^{r} \overline{W}_1^{s+1}){\mathbb E}( W_2^{s}) = {\mathbb E}(W_1^{s} \overline{W}_1^{r+1}){\mathbb E}( W_2^{r})
\en
for non-negative integers $r$ and $s$.  Setting $r=0$, this expresses moments of $W_2$ in terms of the moments of $W_1$.
So the distribution of $W_1$ determines that of $W_2$.  
Assume now the regularity condition {\rm (\ref{p221})}. According to Lemma \ref{beta} we are reduced either to the case 
with $M$ equal frequencies with sum $1$, or to the case where $W_1$ has a beta distribution, and hence so does $W_2$, by
consideration of \re{momeq}. There is nothing more to discuss in the first case, so we assume for the rest of this section
that
\eq
\label{beta12}
\mbox{ each of $W_1$ and $W_2$ has a non-degenerate beta distribution, with possibly different parameters.}
\en
Recall that 
$$
P_1 = W_1 \mbox{ and } P_2 = (1-W_1) W_2.
$$ 
As observed in \cite{jp.isbp}, 
$$
\mbox{the conditional distribution of $(P_3, P_4, \ldots)$ given $P_1$ and $P_2$ depends symmetrically on $P_1$ and $P_2$.}
$$
This can be seen from Kingman's paintbox representation, which implies that conditionally given $\dP_1, \dP_2, \ldots$, as well as $P_1$ and $P_2$, the sequence $(P_3, P_4, \ldots)$ 
is derived by a process of random sampling from the frequencies $(\dP_i)$ with the terms $P_1$ and  $P_2$ deleted.  No matter what $(\dP_i)$ this process depends 
symmetrically on $P_1$ and $P_2$, so the same is true without the extra conditioning on $(\dP_i)$.

Since $P_1+P_2$ is a symmetric function of $P_1$ and $P_2$, and $(W_3,W_4, \ldots$) is a measurable function of $P_1 + P_2$ and $(P_3, P_4, \ldots)$, 
$$
\mbox{ the conditional distribution of $W_3,W_4, \ldots$ given $(P_1, P_2)$ depends symmetrically on $P_1$ and $P_2$.  }
$$
The condition that $W_1$ is independent of $(W_2,W_3, W_4, \ldots)$ implies easily that 
\begin{quote}
$W_1$ is conditionally independent of $(W_3,W_4, \ldots)$ given $W_2$. 
\end{quote}
Otherwise put:
\begin{quote}
$P_1$ is conditionally independent of $(W_3,W_4, \ldots)$ given $P_2/(1-P_1)$, 
\end{quote}
hence by the symmetry discussed above 
\begin{quote}
$P_2$ is conditionally independent of $(W_3, W_4, \ldots)$ given $P_1/(1-P_2)$.
\end{quote}
Let $X:= P_2/(1-P_1)$, $Y:= P_1/(1-P_2)$ and $Z:= (W_3,W_4, \ldots)$. Then we have both
\eq
\label{XZY}
\mbox{ $X$ is conditionally independent of $Z$ given $Y$,}
\en
and
\eq
\label{YZX}
\mbox{ $Y$ is conditionally independent of $Z$ given $X$,}
\en
from which it follows under suitable regularity conditions (see Lemma \ref{JPLemma} below) that
\eq
\label{XYZ}
\mbox{ $(X,Y)$ is independent of $Z$, }
\en
meaning in the present context that
\eq
\label{WWW}
\mbox{ $W_1$, $W_2$ and $(W_3, W_4, \ldots)$ are independent.  }
\en
Lauritzen \cite[Proposition 3.1]{Lauritzen96} shows that \re{XZY} and \re{YZX} imply \re{XYZ}
under the assumption that $(X,Y,Z)$ has a positive and continuous joint density relative to a product measure.
From (\ref{beta12}) and strict positivity of the beta densities on $(0,1)$,
we see that $(X,Y)$ has a strictly positive and continuous density relative to Lebesgue measure on $(0,1)^2$.
We are not in a position to assume that $(X,Y,Z)$ has a density relative to a product measure. However, the passage from
\re{XZY} and \re{YZX} to \re{XYZ} is justified by Lemma \ref{JPLemma} below without need for a trivariate density. 
So we deduce that \re{WWW} holds.
By Lemma \ref{p.prime}, $(W_2, W_3, \ldots)$ is the sequence of residual fractions of an exchangeable partition $\Pi'$, and $W_2$ has a beta density. So either $W_3 = 1$ and we are in the case
\re{alphaneg} with $\couponm = 3$, or $W_3$ has a beta density, and the previous argument applies to show that
\begin{quote}
$W_1$, $W_2$, $W_3$  and $( W_4, W_5, \ldots)$ are independent.  
\end{quote}
Continue by induction to conclude the independence of $W_1, W_2, \ldots W_k$ for all $k$ such that $p(1^k) >0$.
\endpf

\begin{lemma}\label{JPLemma}
Let $X,Y$ and $Z$ denote random variables with values in arbitrary measurable spaces, all defined on a common probability space,
such that {\rm \re{XZY}} and {\rm \re{YZX}} hold.
If the joint distribution of the pair $(X,Y)$ has a strictly positive probability density
relative to some product probability measure, then {\rm \re{XYZ}} holds.
\end{lemma}

\proof
Let $p(X,Y)$ be a version of $\P(Z\in B \mid X,Y)$ for $B$ a measurable set in
the range of $Z$.  
By standard measure theory (e.g.  Kallenberg \cite[6.8]{KallenbergFMP}) the first conditional independence assumption gives
$\P(Z\in B\mid X,Y) = \P(Z \in B \mid X)$ a.s. so that 
\begin{quote}
$p(X,Y)=g(X)$ a.s.  for some measurable function $g$.  
\end{quote}
Similarly from the second conditional independence assumption,
\begin{quote}
$p(X,Y)=h(Y)$ a.s. for some measurable function $h$,
\end{quote}
and we wish to conclude that  
\begin{quote}
$p(X,Y) = c $ a.s. for some constant $c$.
\end{quote}
To complete the argument it suffices to draw this conclusion from the above two assumptions
about a jointly measurable function $p$, with $(X, Y)$ the identity map on 
the product
space of pairs $\mathcal{X} \times \mathcal{Y}$, and the two
almost sure equalities holding with respect to some probability measure $P$ on this space, 
with $P$ having a strictly positive density relative to a product probability  measure $\mu \otimes \nu$.
Fix $u \in (0,1)$, from the previous assumptions it follows that 
\begin{equation}\label{JPLemma1}
\{ p(X,Y) > u\} = \{X  \in A_u\} = \{Y \in C_u\}    ~~~~~ {\rm a.s.}
\end{equation}
for some measurable sets $A_u$, $C_u$, whence 
\begin{equation}\label{JPLemma2}
\{ p(X,Y) > u\} = \{X  \in A_u\} \cap \{Y \in C_u\}    ~~~~~{\rm a.s.},
\end{equation}
where the almost sure equalities hold both with respect to the joint distribution $P$ of $(X,Y)$, and with
respect to a product probability measure $\mu \otimes \nu$ governing $(X,Y)$.  But under $\mu \otimes \nu$ the random variables $X$ and $Y$
are independent. So if $q:= ( \mu \otimes \nu)( p(X,Y) > u)$, then (\ref{JPLemma1}) and 
(\ref{JPLemma2}) imply that $q=q^2$, so $q =0$ or $q=1$.   
Thus $p(X,Y)$ is constant a.s. with respect to $\mu \otimes \nu$, hence also constant with respect to $P$.
\endpf

\section{The deletion property without the regularity condition}
\label{secfails}

Observe that the property required in Theorem \ref{2pchar1} is void if $\pi_n$ happens to be the one-block partition $(n)$.
This readily implies that mixing with 
the trivial one-block partition $\bf 1$ does not destroy the property.
Therefore the $\bf 1$-component may be excluded from the consideration, meaning that it is enough to focus on the case
\eq
\label{nicecase}
\mbox{ $P_1<1$ a.s., or equivalently $\dP_1<1$ a.s., or equivalently $\lim_{n\to\infty} p(n)=0$. }
\en

Suppose then that this condition holds,
but that the first condition in (\ref{p221}) does not hold, so that $p(2,2,1)=0$.
Then 
$$\prob(\dP_2=1-\dP_1>0) +\prob(\dP_1<1,~\dP_2=0)=1.$$
If both terms have positive probability then $\prob(W_2=1\,|\,W_1=0)=0$ but
$\prob(W_2=1\,|\,W_1>0)>0$, so the independence of $W_1$ and $W_2$ fails.
Thus the independence forces either 
$\prob(\dP_2=1-\dP_1>0) =1$ or $\prob(\dP_1<1,~\dP_2=0)=1$. 
The two cases are readily treated:
\begin{itemize}
\item[(i)]
If $\prob(\dP_2=1-\dP_1>0)=1$ then $W_2=1$ a.s. and the independence trivially holds.
This is the case when $\Pi$ has two blocks almost surely.
\item[(ii)]
If  $\prob(\dP_1<1,~\dP_2=0)=1$ and $\prob(\dP_1>0)>0$ then
$\prob(W_2>0\,|\,W_1>0)=0$ but $\prob(W_2>0\,|\,W_1=0)>0$, hence $W_1$ and $W_2$ are not independent.
Therefore   $\prob(\dP_1<1,~\dP_2=0)=1$ and the independence imply $\dP_1=0$ a.s., meaning that $\Pi={\bf 0}$.
\end{itemize}

\par We conclude that the most general exchangeable partition $\Pi$ which has the property 
in Theorem \ref{Kchar1} is a two-component mixture, in which the first component 
is either a partition from the extended two-parameter family, or a two-block partition as in (i) above, or $\bf 0$,
and the second component
is the trivial partition $\bf 1$.

\section{Regeneration and $\tau$-deletion }

In this section we partly survey and partly extend 
the results from 
\cite{RegenComp,rps04}  concerning characterizations of $(\alpha,\theta)$ partitions
by regeneration properties.
As in Kingman's study of  the regenerative processes \cite{KingmanReg},
subordinators (increasing L{\'e}vy processes) appear  naturally in our framework of
{\it multiplicative} regenerative phenomena.  
Following \cite{rps04},
we call a partition structure $(\pi_n)$ {\em regenerative} if 
\begin{quote}
for each $n$ it is possible to delete a randomly chosen part of $\pi_n$ in such a way that for each $0 < m < n$, given the deleted part is of size $m$, the remaining parts form a partition of $n-m$ with the same distribution as $\pi_{n-m}$.
\end{quote}
In terms of an exchangeable partition $\Pi = (\Pi_n)$ of $\mathbb{N}$, the
associated partition structure $(\pi_n)$ is regenerative  if and only if 
\begin{quote}
for each
 $n$ it is possible to select a random block $\B_{nJ_n}$ of $\Pi_n$ in such a way that for each $0 < m < n$, conditionally given that $|\B_{nJ_n}| = m$ the partition $\Pi_n \setminus \B_{nJ_n}$ of $[n-m]$ is distributed according to the unconditional distribution of $\Pi_{n-m}$:
\end{quote}
\eq
\label{regendef}
\mbox{  ($\Pi_n\setminus B_{nJ_n}$ given $|B_{nJ_n}|=m)$ } \ed \Pi_{n-m}
\en
where $\Pi_n \setminus \B_{nJ_n}$ is defined as in the discussion preceding Theorem \ref{Kchar1}.  Moreover, there is no loss of generality in supposing
further that the conditional distribution of $J_n$ given $\Pi_n$ is of the form
\eq
\label{jnjb}
\P(J_n = j \giv \Pi_n = \{B_1, \ldots, B_k\} ) = d(|B_1|,\ldots, |B_k|;j )
\en
for some {\em symmetric deletion kernel} $d$, meaning a non-negative 
function of a composition $\lambda$ of $n$ and $1 \le j \le k_\lambda$
such that
\eq
\label{dnu}
d(\lambda_1,\lambda_2, \ldots, \lambda_k;j ) 
= d(\lambda_{\sigma(1)},\lambda_{\sigma(2)}, \ldots, \lambda_{\sigma(k)};1 )
\en
for every permutation $\sigma$ of $[k]$ with $\sigma(1) = j$.
To determine a symmetric deletion kernel, is suffices to specify 
$d(\lambda;1 )$, 
which is the conditional probability, given blocks of sizes $\lambda_1,\lambda_2, \ldots, \lambda_k$,
 of picking the first of these blocks. This is a non-negative symmetric function of
$(\lambda_2, \ldots, \lambda_k)$, subject to the further constraint that its extension 
to arguments $j \ne 1$ via \re{dnu} satisfies
$$
\sum_{j = 1}^{k_\lambda} d(\lambda; j ) = 1
$$
for every composition $\lambda$ of $n$.
The regeneration condition can now be reformulated 
in terms of the EPPF $p$ of $\Pi$ in a manner similar to \re{factor}:

\begin{lemma}
An exchangeable random partition $\Pi$ with EPPF $p$ is regenerative 
if and only if there exists a symmetric deletion kernel $d$ such that
\eq
\label{dnu1}
p( \lambda )
d( \lambda ; 1)
= q(n,\lambda_1)
{n \choose \lambda_1}^{-1} p(\lambda_2, \ldots, \lambda_k)
\en
for every composition $\lambda$ of $n$ into at least two parts and some non-negative function $q$.
Then 
\eq
\label{decmat}
q(n,m) = \P(|\B_{n,J_n}| = m  )     ~~~~(m \in [n])
\en
for $J_n$ as in {\rm \re{jnjb}}.
\end{lemma}
\proof
Formula \re{dnu1} offers two different ways of computing the probability
of the event that $\Pi_n = \{B_1, \ldots, B_k \}$ and $J_n = 1$ for
an arbitrary partition $\{B_1, \ldots, B_k \}$ of $[n]$ with $|B_i| = \lambda_i$
for $i \in [k]$: on the left side, by definition of the symmetric deletion
kernel, and on the right side by conditioning on the event 
$\B_{n,J_n} = B_1$ and appealing to the regeneration property and 
exchangeability.
\endpf

Consider now the question of whether an $(\alpha,\theta)$ partition with EPPF 
$p = p_{\alpha,\theta}$ as in \re{EPPF2Par} is regenerative with respect to some
deletion kernel. By the previous lemma and cancellation of common factors, the question is whether there exists
a symmetric deletion kernel $d(\lambda;j)$ such that the value of
\eq
\label{f1}
q(n,\lambda_1)= d(\lambda; 1) {n \choose \lambda_1} \frac{
(1 - \alpha)_{\lambda_1 - 1}
 ( \theta + (k-1) \alpha ) 
}{ ( \theta + n - \lambda_1 )_{\lambda_1} }
\en
is the same for all compositions $\lambda$ of $n$ with $k$ parts and a prescribed 
value of $\lambda_1$. But it is easily checked that  the formula
\eq
\label{f2}
d_{\alpha,\theta}(\lambda ; j) = 
\frac{ \theta \lambda_j + \alpha (n - \lambda_j) }{ n (\theta + \alpha (k-1) ) }
\en
provides just such a symmetric deletion kernel. Note that the kernel depends
on $(\alpha,\theta)$ only through the ratio $\tau:= \alpha/(\alpha + \theta)$,
and that the kernel is non-negative for all compositions $\lambda$ only if both $\alpha$ and $\theta$
are non-negative.

To provide a more general context for this and later discussions, let $( x_1, \ldots, x_k)$ be a fixed sequence of 
positive numbers with sum $s = \sum_{j=1}^k x_j$.  For a fixed parameter $\tau\in [0,1]$, define a random variable $T$ with values in $[k]$ by 
\begin{equation}\label{T-dist} 
\prob\left(T=j \mid (x_1, \ldots, x_k)\right)=
{(1-\tau)x_j + \tau (s-x_j)
\over s(1-\tau+\tau (k-1) )},
\end{equation}
The random variable $x_{T}$ is called a {\it $\tau$-biased pick} from $x_1, \ldots, x_k$.
The law of $x_{T}$ does not depend on the order of the sequence $(x_1, \ldots, x_k)$, and there is also a scaling invariance: $s^{-1}x_T$ is a $\tau$-biased pick from $(s^{-1}x_1, \ldots, s^{-1}x_k)$.  
Note that a $0$-biased pick is a size-biased pick from $(x_1, \ldots, x_k)$, choosing any particular element with probability proportional to its size.  A $1/2$-biased pick is 
a uniform random choice from the list, as (\ref{T-dist}) then equals $1/k$ for all $j$.  And 
a $1$-biased pick may be called a co-size biased pick, as it 
chooses $j$ with probability proportional to its co-size $s-x_j$.

\par These definitions are now applied to the sequence of block sizes $x_j$ of the restriction to $[n]$  of an an exchangeable partition 
$\Pi$ of $\mathbb{N}$.
We denote by $T_n$ a random variable whose conditional distribution 
given $\Pi_n$ with $k$ blocks and $|\B_{nj}| = x_j$ for $j \in [k]$ is
defined by (\ref{T-dist}), and denote by $\B_{nT_n}$ the $\tau$-biased pick from
the sequence of blocks of $\Pi_n$. We call $\Pi$ {\it $\tau$-regenerative} if  
$\Pi_n$ is regenerative with respect to deletion of the $\tau$-biased pick $\B_{nT_n}$.  

\begin{theorem}{\rm \cite{ RegenComp, rps04}} \label{T-char}
For each $\tau \in [0,1]$, 
apart from the constant partitions $\mathbf{0}$ and $\mathbf{1}$, 
the only exchangeable partitions of $\mathbb{N}$
that are $\tau$-regenerative are the members of the two parameter family with parameters in the range 
$$\{(\alpha,\theta)\in [0,1]\times [0,\infty]: ~ \alpha/(\alpha+\theta)=\tau\}.$$
Explicitly, the distribution of the $\tau$-biased pick for such $(\alpha,\theta)$ partitions of $[n]$ is
\eq
\label{2paramdec}
\prob(|\B_{nT_n}|=m)={n\choose m}{(1-\alpha)_{m-1}\over (\theta+n-m)_{m}}
{(n-m)\alpha+ m\theta\over n},~~~~~~m\in[n].
\en
\end{theorem}
\proof
The preceding discussion around \re{f1} and \re{f2} shows that
members of the two parameter family with parameters in the indicated
range  are  $\tau$-regenerative, and gives the formula \re{2paramdec} for the decrement matrix.
See \cite{rps04} for the proof that these are the only non-degenerate 
exchangeable partitions of $\mathbb{N}$ that are $\tau$-regenerative.
\endpf

In particular, each $(\alpha,\alpha)$ partition is $1/2$-regenerative, 
meaning regenerative with respect to deletion of a
block chosen uniformly at random.  
The constant partitions $\mathbf{0}$ and $\mathbf{1}$ 
are obviously $\tau$ regenerative for every $\tau \in [0,1]$.  This is consistent with the characterization above because the $(1,\theta)$ partition is the $\mathbf{0}$ partition for every $\theta \geq 0$, and because the partition $\mathbf{1}$ can be reached as a limit of $(\alpha, \theta)$ partitions as $\alpha, \theta \downarrow 0$ 
with $\alpha(\alpha+\theta)^{-1}$ held fixed.

\paragraph{Multiplicative regeneration} 
 
By Corollary \ref{2pcharcor}, if $(P_i)$ is the sequence of limit frequencies for a $(0,\theta)$ partition for some $\theta>0$ and if the first limit frequency $P_1$ is deleted and the other frequencies renormalized to sum to 1, then the resulting sequence $(Q_j)$ is independent of $P_1$ and has the same 
distribution as $(P_j)$.  Because $P_1$ is a size-biased pick from the sequence $(P_i)$, this regenerative property of the frequencies $(P_i)$ can be seen as an analogue of the $0$-regeneration property of the $(0, \theta)$ partitions.  

If $(P_i)$ is instead the sequence of limit frequencies of an $(\alpha, \theta)$ partition $\Pi$ for parameters satisfying $0<\alpha<1, ~\alpha/(\alpha+\theta)=\tau$,
a question arises: does the regenerative property of $\Pi_n$ with respect to a $\tau$-biased pick have an analogue in terms of a $\tau$-biased pick from the frequencies $(P_i)$?
This cannot be answered  straightforwardly as in the $\tau=0$ case,
 because when $\tau>0$ the formula
(\ref{T-dist}) defines a proper probability distribution 
only for series $(x_j)$ with some finite number $k$
of positive terms. 
For instance, in the case $\tau=1/2$ there is no such analogue of (\ref{T-dist})
as `uniform random choice' from infinitely many terms.

\par Still, Ewens' case provides a clue if we turn to a {\it bulk deletion}. 
Let $P_J$ be a size-biased pick from the frequencies $(P_j)$, as 
defined by (\ref{sbpick}), and let $(Q_j)$ be a sequence obtained from $(P_j)$ by deleting all $P_1,\ldots,P_J$ and renormalizing.
Then $(Q_j)$ is independent of $P_1,\dots,P_J$, and $(Q_j)\ed(P_j)$. The latter assertion follows from the i.i.d. property 
of the residual fractions and by noting that
(\ref{sbpick}) is identical with 
$$\prob(J=j\,|\, (W_i, i\in \Nat))= W_j\prod_{i=1}^{j-1}(1-W_i).$$
A similar bulk deletion property holds for partitions in the Ewens' family, in the form:
$$
\mbox{  ($\Pi_n\setminus (\B_{n1}\cup\dots\cup\B_{nJ_n})$ given $|\B_1\cup\dots\cup\B_{nJ_n}|=m)$ } \ed \Pi_{n-m}
$$
for all $1 \leq m \leq n$, where $B_{nJ_n}$ is a size-biased pick from the
blocks.

To make the ansatz of bulk deletion work for $\tau\neq 0$ it is necessary to arrange the frequencies in a
more complex manner. 
To start with, we modify the paintbox construction.
Let ${\cal U}\subset [0,1]$ be a random open set canonically represented  as the union of its disjoint open
component intervals. 
We suppose that the Lebesgue measure of $\cal U$, equal to the sum of lengths of the components,
is $1$ almost surely.  We associate with $\cal U$ an exchangeable partition 
$\Pi$ 
exactly as in Kingman's representation
in Theorem \ref{Kpaintbox}. 
For each component interval $G \subset \cal U$ there is an index 
$i_G:= 
\min\{n: U_n \in G \}$ that is the minimal index 
of a sequence $(U_i)$ of iid uniform[0,1] points hitting the interval, 
and for all $j$, $P_{j}$ is the  length of the $j$th component interval when the intervals 
are listed in order of increasing minimal indices.  So $(P_{j})$ is a size-biased permutation of 
the lengths of interval components of $\cal U$.

Let $\triangleleft$ be the linear order on $\Nat$ induced by the 
{\it interval order} of the components of $\cal U$,  so $j\triangleleft k$ iff 
the interval of length $P_j$, which is the home interval of the $j$th block $\B_j$ to appear in the process of uniform random
sampling of intervals, lies to the left of the interval of length $P_k$ associated with block $\B_k$.
A convergence argument shows that 
$\cal U$ is uniquely determined by $(P_j)$ and $\triangleleft$.
In loose terms,   $\cal U$ is an arrangement of a sequence of tiles of sizes $P_j$ in the order on indices $j$
prescribed by $\triangleleft$, and this arrangement is constructable by sequentially
placing the tile $j$ in the position prescribed by the order $\triangleleft$ restricted to $[j]$.

For $x\in [0,1)$ let $(a_x,b_x)\subset{\cal U}$ be the component interval
containing $x$. Define $\mathcal{V}_x$ as the open set
obtained by deleting
the bulk of component intervals to the left of $b_x$, then linearly rescaling the remaining set 
${\mathcal U}\cap [b_x,1]$ to $[0,1]$. 
We say that $\mathcal U$ is {\it multiplicatively regenerative} if for each $x\in [0,1)$,
${\mathcal{V}_x}$ is independent of $\mathcal{U} \cap [0,b_x]$ and $\mathcal{V}_x \ed {\cal U}$.

\par An ordered version of the paintbox correspondence yields:  

\begin{theorem}\label{PartPa}{\rm \cite{ RegenComp, rps04}} An exchangeable partition $\Pi$ is regenerative if and only if
it has a paintbox  representation in terms of some multiplicatively regenerative set $\cal U$.
The deletion operation is then defined by classifying $n$ independent uniform points
 from $[0,1]$ according to the intervals of ${\cal U}$ into which they fall, and deleting 
the block of points in the leftmost occupied interval.
\end{theorem}

A property of the frequencies $(P_j)$ of an exchangeable regenerative partition $\Pi$ of $\Nat$
now emerges: there exists a strict total order $\triangleleft$ on $\Nat$, which is a random order, which has some joint distribution with
$(P_j)$ such that arranging the intervals of sizes $(P_j)$ in order $\triangleleft$ yields a multiplicatively regenerative set ${\cal U}$.
Equivalently, there exists a multiplicatively regenerative set ${\cal U}$ that induces a partition
with frequencies $(P_j)$ and an associated order $\triangleleft$.  This set ${\cal U}$ is then necessarily unique 
in distribution as a random element 
of the space of open subsets 
of $[0,1]$ equipped with the Hausdorff metric \cite{RegenComp} on the complementary closed subsets. 
A subtle point here is that the joint distribution of $(P_j)$ and $\triangleleft$ is not unique,
and neither is the joint distribution of $(P_j)$ and ${\cal U}$, unless further conditions are imposed.
For instance, one way to generate $\triangleleft$ is to suppose that the $(P_j)$ are generated by a process of uniform random
from ${\cal U}$. But for a $(0,\theta)$ partition, we know that another way is to construct ${\cal U}$ from
$(P_j)$ by simply placing the intervals in deterministic order $P_1, P_2, \ldots$ from left to right.
In the construction by uniform random sampling from $\cal U$ the interval of length $P_1$ discovered by the first  sample point
need not be the leftmost, and need not lie to the left of the second discovered interval $P_2$.

In \cite{RegenComp} we showed that the multiplicative regeneration of  $\mathcal{U}$ follows from 
an apparently weaker property:
if $(a_U, b_U)$ is the component interval of $\cal U$ containing an uniform[0,1] 
sample $U$ independent of $\mathcal{U}$, and if $\mathcal{V}$ is defined as 
the open set obtained by deleting the component intervals to the left of $b_U$ 
and linearly rescaling the remaining set $\mathcal{U}\cap [b_U,1]$ to [0,1], then given $b_U<1$, 
$\mathcal{V}$ is independent of $b_U$ 
(hence, as we proved, independent of $\mathcal{U}\cap [0,b_U]$ too!)
and has distribution equal 
to the unconditional distribution of $\mathcal{U}$.  
This independence is the desired analogue for more general regenerative partitions of the 
bulk-deletion property of Ewens' partitions.

\par The fundamental representation of multiplicatively regenerative sets
involves a random process $F_t$ known in statistics as a neutral-to-the right distribution function.

\begin{theorem}\label{m-reg-rep} {\rm \cite{RegenComp}}
A random open set $\cal U$ of Lebesgue measure $1$ is  multiplicatively regenerative 
if and only if there exists a drift-free subordinator $S=(S_t,t\geq 0)$ with $S_0=0$ such that $\cal U$
is the complement to the closed range of the process $F_t=1-\exp(-S_t),~t\geq 0$.
The L{\'e}vy measure of $S$ is determined uniquely up to a positive factor.
\end{theorem}
According to Theorems \ref{PartPa} and \ref{m-reg-rep}, regenerative partition structures with 
proper frequencies are parameterised by a measure $\nutil({\rm d}u)$ on $(0,1]$ with 
finite first moment, which is the image 
via the transformation from $s$ to $1 - \exp(-s)$
of the L{\'e}vy measure $\nu({\rm d}s)$ on $(0,\infty]$
associated with the subordinator $S$.
The Laplace exponent $\Phi$ of the subordinator, 
defined 
by the L\'evy-Khintchine formula
$$
\E [ \exp ( - a S_t ) ] =  \exp [ - t \Phi(a) ],   ~~~~~~~~a \ge 0
$$
determines the L\'evy measure $\nu({\rm d}s)$ on $(0,\infty]$
and its image $\nutil({\rm d}u)$  on $(0,1]$ via the formulae
$$
\Phi(a)= 
\int_{(0,\infty]} (1 - e^{-a x}) \nu(dx) 
= \int_{]0,1]}(1-(1-x)^a)\nutil({\rm d}x).
$$
As shown in \cite{RegenComp}, the decrement matrix $q$ of the regenerative partition structure,
as in \re{decmat}, is then
$$
q(n,m)={\Phi(n,m)\over \Phi(n)}\,,\qquad 1\leq m\leq n\,,~n=1,2,\ldots
$$
where
$$
\Phi(n,m)= {n\choose m}\int_{]0,1]} x^m(1-x)^{n-m}\nutil({\rm d}x)\,.
$$ 
Uniqueness of the parameterisation is achieved by a normalisation condition, such as $\Phi(1)=1$.

In \cite{RegenComp} the subordinator $S^{\alpha,\theta}$
which produces $\cal U$ 
as in Theorem \ref{m-reg-rep} for the $(\alpha, \theta)$ partition was identified by the following formula
for the right tail  of
its L{\'e}vy measure: 
\begin{equation}\label{LM}
\nu^{}(x,\infty]= (1-e^{-x })^{-\alpha} e^{-x\theta}, ~~~x>0.
\end{equation}

The subordinator $S^{(0,\theta)}$
is a compound Poisson process whose jumps are exponentially distributed
with rate $\theta$. 
For $\theta=0$ the L{\'e}vy measure has a unit mass at $\infty$, so 
the subordinator $S^{(\alpha,0)}$ is killed at unit rate. 
The $S^{(\alpha,\alpha)}$ 
subordinator belongs to the class of Lamperti-stable processes recently studied in \cite{Lstable}.
For positive parameters the subordinator $S^{(\alpha,\theta)}$  can be constructed from 
the $(0,\theta)$ and $(\alpha,0)$ cases, as follows.
First split ${\mathbb R}_+$ by the range of $S^{(0,\theta)}$, that is at points $E_1<E_2<\dots$ 
of a Poisson process with rate $\theta$. Then run an independent copy of $S^{(\alpha,0)}$ 
up to the moment the process crosses $E_1$ at some random time,
say $t_1$. The level-overshooting value is neglected and the process is stopped. At the same time
$t_1$ a new independent copy of $S^{(\alpha,0)}$ 
is started at value $E_1$ and run until crossing $E_2$ at some random time $t_2$, and so on. 

In terms of $F_t=1-\exp(-S_t)$, the range of the process in the $(0,\theta)$ case is a
stick-breaking set $\{1-\prod_{i=1}^{j-1} (1-V_i), i=0,1,\ldots\}$ 
with i.i.d. beta$(1,\theta)$ factors $V_i$. 
In the case $(\alpha, 0)$ 
the range of $(F_t)$ is 
the intersection of $[0,1]$ with the $\alpha$-stable set
(the range of $\alpha$-stable subordinator). 
In other cases $\cal U$ is constructable as a cross-breed of the cases $(\theta,0)$ and 
$(0,\alpha)$:
first $[0,1]$ is partitioned in subintervals by the beta$(1,\theta)$ stick-breaking,
then each subinterval $(a,b)$ of this partition is further split by independent copy of the multiplicatively regenerative
$(\alpha,0)$ set, shifted to start at $a$ and truncated at $b$.

\paragraph{Constructing the order} 
Following \cite{rps04, PitmanWinkel}, we shall describe 
an arrangement
which allows us to pass from $(\alpha, \theta)$ frequencies $(P_j)$ to the multiplicatively regenerative set
associated with the subordinator $S^{(\alpha,\alpha)}$. 
The connection between size-biased permutation with $\tau$-deletion (Lemma \ref{Leem}) is new.

A linear order $\triangleleft$ on $\Nat$ is conveniently described by a sequence of the initial ranks 
$(\rho_j)\in [1]\times[2]\times\cdots$, with $\rho_j=i$ if and only if $j$ is ranked $i$th smallest 
in the order $\triangleleft$ among the integers $1,\dots,j$.
For instance, the initial ranks $1,2,1,3\dots$ appear when $3\triangleleft 1\triangleleft 4\triangleleft 2$.

For $\xi\in[0,\infty]$  define a random order $\triangleleft_\xi$ on $\Nat$ by 
assuming that the initial ranks $\rho_k, k\in \Nat,$ are independent, with distribution 
$$\prob(\rho_k=j)=
  {1\over k+\xi-1}{\tt 1}(0 < j<k)
+
{\xi\over k+\xi-1}{\tt 1}(j=k) ~~~~~~~~~~~,k>1.
$$
The edge cases $\xi=0,\infty$ are defined by continuity. The order $\triangleleft_1$
is a `uniformly random order', in the sense that restricting to $[n]$ we have all $n!$
orders equally likely, for every $n$. The order $\triangleleft_\infty$ coincides
with the standard order $<$ almost surely.
For every permutation $i_1,\ldots,i_n$ of $[n]$, we have
$$\prob(i_1\triangleleft_\xi\dots\triangleleft_\xi i_n)={\xi^r\over \xi(\xi+1)\dots(\xi+n-1)}$$
where $r$ is the number of upper records in the permutation.
See \cite{Coherent} for this and more general permutations with tilted record statistics.

\begin{theorem}\label{newthm} 
{\rm \cite[Corollary 7]{PitmanWinkel} }
For $0\leq \alpha<1, \theta\geq0$ the arrangement 
of  $GEM(\alpha,\theta)$ frequencies $(P_j)$ represented as open intervals in an independent random order $\triangleleft_{\theta/\alpha}$ is a
multiplicatively regenerative open set ${\cal U}\subset[0,1]$,
where $\cal U$ is representable
as the complement of the closed range of
the process $F_t=1-\exp(-S_t), t\geq 0,$ for $S$ the subordinator with L{\'e}vy measure {\rm (\ref{LM})}.
\end{theorem}


This result was presented without proof as \cite[Corollary 7]{PitmanWinkel},
in a context where the regenerative ordering of frequencies was motivated by 
an application to a tree growth process. 
Here we offer a proof which exposes the combinatorial structure of the
composition of size-biased permutation and  a $\triangleleft_{\theta/\alpha}$ ordering of frequencies.

For a sequence of positive reals $(x_1,\dots,x_k)$, define the {\it $\tau$-biased permutation} 
of this sequence, denoted ${\rm perm}_\tau(x_1,\ldots,x_k)$,
by iterating  a single $\tau$-biased pick, as follows.  A number $x_{T}$ is chosen  from $x_1,\dots,x_k$
 without replacement,
with $T$ distributed on $[k]$ according to (\ref{T-dist}), and $x_T$ is placed in position $1$. 
Then the next number is chosen from $k-1$ remaining numbers
using again the rule of $\tau$-biased pick, 
and placed in position 2, etc. 

The instance ${\rm perm}_0$ is the size-biased permutation, which
is defined more widely for
 finite or infinite summable sequences $(x_1,x_2,\ldots)$, and shuffles them 
in the same way as it shuffles $(s^{-1}x_1, s^{-1}x_2,\dots)$ where $s=\sum_{j} x_j$.
Denote by $\triangleleft_\xi(x_1,\dots,x_k)$ the arrangement of $x_1,\ldots,x_k$ in succession according to the
 $\triangleleft_\xi$-order on $[k]$. 
 
\begin{lemma}\label{Leem} For $\xi=(1 - \tau)/\tau$ there is the compositional formula
\eq\label{coin}
{\rm perm}_\tau(x_1,\ldots,x_k)\ed {\triangleleft}_\xi({\rm perm}_0(x_1,\ldots,x_k)),
\en
where on the right-hand side 
$\triangleleft_\xi$ 
and 
${\rm perm}_0$ 
are independent.
\end{lemma}
\proof 
On each side of this identity, the distribution of the random permutation remains the same if the
sequence $x_1, \ldots, x_k$ is permuted. So it suffices to check that each scheme 
returns the identity permutation with the same probability. If on the right hand side we set 
$$
{\rm perm}_0(x_1,\ldots,x_k) = ( x_{\sigma(1)},  \ldots, x_{\sigma(k)})
$$
then the right hand scheme generates the identity permutation with probability
\eq
\label{ExiR}
\frac{ \E \xi^{R}  } {\xi ( \xi +1 ) \cdots (\xi + k - 1 )}
\en
where $R$ is the number of upper records in the sequence of ranks which generated $\sigma^{-1}$, which equals the number of upper records in $\sigma$. Now
$R = \sum_{j = 1}^k X_j$ where $X_j$ is the indicator of the event $A_j$ that $j$ is an upper record level for $\sigma$, meaning that
there is some $1 \le i \le n$ such that 
$$
\mbox{ $\sigma(i') <  j $ for all $i' < i$ and $\sigma(i) = j$.}
$$
Equivalently, $A_j$ is the event that 
$$
\mbox{ $\sigma^{-1}(j) < \sigma^{-1}(\ell)$ for each $j < \ell \le k$. }
$$
Or again, assuming for simplicity
that the $x_i$ are all distinct, which involves no loss of generality, because the probability in question depends continuously on $(x_1, \ldots, x_k)$,
$A_j$ is the event that $x_j$ precedes $x_\ell$ in the permutation $( x_{\sigma(1)},  \ldots, x_{\sigma(k)})$
for each $j < \ell \le k$.
Now it is easily shown that $( x_{\sigma(1)},  \ldots, x_{\sigma(k)})$ with $x_1$ deleted is a size-biased permutation of $(x_2, \ldots, x_k)$, and
that the same is true conditionally given $A_1$. It follows by induction that the events $A_j$ are mutually independent, with
$$
\P(A_j) = x_j/(x_j + \cdots + x_k)            \mbox{ for } 1 \le j \le k.
$$
This allows the probability in \re{ExiR} to be evaluated as
$$
\prod_{j = 1}^k  \frac{(\xi x_j + x_{j+1} + \cdots +  x_k ) } { ( x_j + x_{j+1} + \cdots + x_k)(\xi + j -1 ) }
$$
This is evidently the probability that ${\rm perm}_\tau(x_1,\ldots,x_k)$ generates the identity permutation,
and the conclusion follows.
\endpf

\par The $\tau$-biased arrangement cannot be defined for infinite positive summable sequence
$(x_1,x_2,\dots)$, since the `$k=\infty$' instance of (\ref{T-dist}) is not a proper distribution for $\tau\neq 0$.
But the right-hand side of (\ref{coin}) is well-defined
as arrangement of $x_1,x_2,\dots$ in some total order,
hence   
the composition ${\triangleleft}_\xi\circ{\rm perm}_0$ is the natural extension of 
the $\tau$-biased arrangement to infinite series.
\vskip0.3cm

\noindent
{\it Proof of Theorem \ref{newthm}.} We represent a finite or infinite positive sequence $(x_j)$ whose sum is  $1$
as an open subset of $[0,1]$ composed of contiguous intervals of sizes $x_j$. The space of open subsets of $[0,1]$
is endowed with the Hausdorff distance on the complementary compact sets. This topology is weaker than the product topology 
on positive series summable to $1$. 
The limits below are understood as $n\to\infty$.

\par We know by a version of
Kingman's correspondence \cite{jp.epe} that $(|\B_{nj}|/n, j\geq 1)\to (P_j)$ a.s. in the product topology. 
This readily implies $\triangleleft_\xi (|\B_{nj}|/n, j\geq 1)\to \triangleleft_\xi(P_j)$ a.s.
in the Hausdorff topology, by looking at the $M$ first terms for $M$ such that these terms sum to at least  $1-\epsilon$
with probability at least $1 -\epsilon$, then sending $\epsilon \to 0$ and $M\to\infty$.
In \cite{RegenComp} we showed that ${\rm perm}_\tau(|\B_{nj}|, j\geq 1)\to {\cal U}$ a.s. in the Hausdorff topology.
(Here the definition of the ${\rm perm}_\tau$ is coupled with $(|\B_{nj}|, j\geq 1)$ by putting these blocks in the order
determined by uniform sampling from ${\cal U}$).
The missing link is provided by Lemma \ref{Leem}, from which  we obtain
$${\rm perm}_\tau(|\B_{nj}|, j\geq 1)\ed \triangleleft_\xi(|\B_{nj}|,j\geq 1),$$
with the $\tau$-biased permutation ${\rm perm}_\tau$ applied to the {\it finite} sequence of positive 
block-sizes $(|\B_{nj}|, j\geq 1)$. 
Putting things together we conclude that 
$\triangleleft_\xi(P_j, j\ge 1)\ed{\cal U}$.
\endpf

In three special cases, already identified in the previous work \cite{rps04}, 
the arrangement of PD$(\alpha,\theta)$ (or GEM$(\alpha,\theta)$) frequencies 
in a multiplicatively regenerative set has a simpler description:
in the $(0,\theta)$ case the frequencies are placed in the size-biased order;
in the $(\alpha,\alpha)$ case the frequencies are `uniformly randomly shuffled';
and in the $(\alpha,0)$ case a size-biased pick is placed 
contiguously to 1, while the other frequencies are `uniformly randomly shuffled'.
The latter is an infinite analogue of the co-size biased arrangement ${\rm perm}_1$.

\par We refer to \cite{hpw07, PitmanWinkel} for further recent developments related to 
 ordered $(\alpha,\theta)$ partitions and their regenerative properties.

\end{document}